\documentclass[11pt]{amsart}

\allowdisplaybreaks[4]
\newtheorem{proposition}{Proposition}[section]

\newtheorem{corollary}[proposition]{Corollary}
\newtheorem{theorem}[proposition]{Theorem}

\theoremstyle{definition}
\newtheorem{definition}[proposition]{Definition}
\newtheorem{example}[proposition]{Example}

\newtheorem{remark}[proposition]{Remark}
\newtheorem{remarks}[proposition]{Remarks}
\newtheorem{problem}[proposition]{Problem}

\newcommand{\thlabel}[1]{\label{th:#1}}
\newcommand{\thref}[1]{Theorem~\ref{th:#1}}
\newcommand{\selabel}[1]{\label{se:#1}}
\newcommand{\seref}[1]{Section~\ref{se:#1}}

\newcommand{\prlabel}[1]{\label{pr:#1}}
\newcommand{\prref}[1]{Proposition~\ref{pr:#1}}
\newcommand{\colabel}[1]{\label{co:#1}}
\newcommand{\coref}[1]{Corollary~\ref{co:#1}}
\newcommand{\relabel}[1]{\label{re:#1}}
\newcommand{\reref}[1]{Remark~\ref{re:#1}}
\newcommand{\exlabel}[1]{\label{ex:#1}}

\newcommand{\delabel}[1]{\label{de:#1}}
\newcommand{\deref}[1]{Definition~\ref{de:#1}}
\newcommand{\eqlabel}[1]{\label{eq:#1}}
\newcommand{\equref}[1]{(\ref{eq:#1})}

\newcommand{\Hom}{{\rm Hom}}
\newcommand{\End}{{\rm End}}

\newcommand{\im}{{\rm Im}\,}

\def\lan{\langle}
\def\ran{\rangle}
\def\ot{\otimes}

\newcommand{\Aa}{\mathcal{A}}

\newcommand{\Cc}{\mathcal{C}}

\newcommand{\Mm}{\mathcal{M}}

\newcommand{\YD}{\mathcal{YD}}
\newcommand{\Ww}{\mathcal{W}}
\newcommand{\Zz}{\mathcal{Z}}

\def\*C{{}^*\hspace*{-1pt}{\Cc}}
\def\text#1{{\rm {\rm #1}}}

\def\ul{\underline}

\def\equal#1{\smash{\mathop{=}\limits^{#1}}}

\usepackage{xypic}

\usepackage{color,amssymb,graphicx}

\begin{document}
\title[The center of the category of bimodules]
{The center of the category of bimodules and descent data for non-commutative rings}

\author{A. L. Agore}
\address{Faculty of Engineering, Vrije Universiteit Brussel, Pleinlaan 2, B-1050 Brussels, Belgium}
\email{ana.agore@vub.ac.be and ana.agore@gmail.com}

\author{S. Caenepeel}
\address{Faculty of Engineering, Vrije Universiteit Brussel, Pleinlaan 2, B-1050 Brussels, Belgium}
\email{scaenepe@vub.ac.be}

\author{G. Militaru}
\address{Faculty of Mathematics and Computer Science, University of Bucharest, Str.
Academiei 14, RO-010014 Bucharest 1, Romania}
\email{gigel.militaru@fmi.unibuc.ro and gigel.militaru@gmail.com}


\thanks{A.L.A. is research fellow ''aspirant'' of FWO-Vlaanderen. S.C. is supported by
FWO project G.0117.10 ``Equivariant Brauer groups and Galois
deformations''. G.M. is supported by the CNCS - UEFISCDI grant no.
88/05.10.2011 ''Hopf algebras and related topics''.}

\begin{abstract}
Let $A$ be an algebra over a commutative ring $k$. We compute the
center of the category of $A$-bimodules. There are six isomorphic
descriptions: the center equals the weak center, and can be
described as categories of noncommutative descent data, comodules
over the Sweedler canonical $A$-coring, Yetter-Drinfeld type
modules or modules with a flat connection from noncommutative
differential geometry. All six isomorphic categories are braided
monoidal categories: in particular, the category of comodules over
the Sweedler canonical $A$-coring $A \ot A$ is braided monoidal.
We provide several applications: for instance, if $A$ is finitely
generated projective over $k$ then the category of left
$\End_k(A)$-modules is braided monoidal and we give an explicit
description of the braiding in terms of the finite dual basis of
$A$. As another application, new families of solutions for the
quantum Yang-Baxter equation are constructed: they are canonical
maps $\Omega$ associated to any right comodule over the Sweedler
canonical coring $A \ot A$ and satisfy the condition $\Omega^3 =
\Omega$. Explicit examples are provided.
\end{abstract}

\maketitle

\section*{Introduction}
A monoidal category can be viewed as a categorical version of a
monoid. The appropriate generalization of the center of a monoid
is given by the centre construction, which was introduced
independently by Drinfeld (unpublished), Joyal and Street
\cite{JS2} and Majid \cite{m2}. A key result in the classical
theory is the following: the center of the category of
representations of a Hopf algebra $H$ is isomorphic to the
category of Yetter-Drinfeld modules over $H$ \cite{K}. Moreover,
if the Hopf algebra $H$ is finite diminesional, then the category
of Yetter-Drinfeld modules is isomorphic to the category of
representations over the Drinfeld double $D(H)$. Since the center
is a braided monoidal category, it follows that the Drinfeld
double is a quasitriangular Hopf algebra.

Let $A$ be an algebra over a commutative ring $k$. In this note,
we study the center of the category ${}_A\Mm_A$ of $A$-bimodules,
and relate it to some classical concepts. We introduce $A\ot
A^{\rm op}$-Yetter-Drinfeld modules (\deref{3.0}), and show that
the weak center of ${}_A\Mm_A$ is isomorphic to the category of
$A\ot A^{\rm op}$-Yetter-Drinfeld modules (\prref{3.3}). We give
other descriptions: the weak center is equal to the center
(\prref{3.5}) and is isomorphic to the category $\Mm^{A\ot A}$ of
comodules over the Sweedler canonical coring $A\ot A$
(\prref{3.2}). Moreover it was proved in \cite[Theorem 5]{tb06}
that the category $\Mm^{A\ot A}$ is isomorphic to the category of
right $A$-modules with a flat connection as defined in
noncommutative differential geometry. Thus, the category of right
$A$-modules with a flat connection is also equal to the center. We
introduce a category of descent data $\ul{\rm Desc}(A/k)$,
generalizing the descent data introduced in \cite{KO} from $A$
commutative to $A$ non-commutative, and this category is also
isomorphic to the center. The first main result of this paper is
summarized in \thref{maincentru} which provide six isomorphic
descriptions for the center of the category of $A$-bimodules. All
six isomorphic categories are braided monoidal categories. In
particular, the category of comodules over the Sweedler canonical
$A$-coring $A\ot A$ is a braided monoidal category and hence one
can perform most of the constructions that are performed for
differentiable manifolds. For instance, connections in bimodules
try to mimic linear connections in geometry and are useful in
capturing Riemannian aspects (see \cite{BrzezinskiWisbauer},
\cite{tb11} for more detalis). If $A$ is faithfully flat as a
$k$-module, all these categories are equivalent to the category of
$k$-modules, by classical descent theory. In the case where $A$ is
finitely generated and projective, the category $\Mm^{A\ot A}$ is
isomorphic to the category of left modules over $\End_k(A)$, in
fact, one may view $\End_k(A)$ as the Drinfeld double of the
enveloping algebra $A^e = A\ot A^{\rm op}$. Thus, the category of
left $\End_k(A)$-modules is braided monoidal, and we give an
explicit description of the monoidal structure and the braiding in
terms of the finite dual basis of $A$. If we apply this to the
case where $A=k^n$, then we find that the category of left modules
over the matrix ring $M_n(k)$ is braided monoidal. We give an
explicit description of the tensor product and the braiding.

The second major application of the above results is the fact that
they lead to constructing new and interesting family of solutions
for the quantum Yang-Baxter equation. If $V$ is a right comodule
over the Sweedler canonical coring $A \ot A$, then the canonical
map $\Omega : V\ot V \to V\ot V$, $\Omega (v \ot w) = w_{[0]} \ot
v_{[0]} w_{[1]} v_{[1]}$ is a solution of the quantum Yang-Baxter
equation and $\Omega^3 = \Omega$ in the endomorphism algebra $\End
(V \ot V)$ (\thref{3.21}). Several examples are provided.

\section{Preliminary Results}
\subsection{Braided monoidal categories and the center construction}\selabel{0.1}
A monoidal category $\Cc=(\Cc,\ot, I,a,l,r)$ consists of a
category $\Cc$, a functor $\ot:\ \Cc\times\Cc\to \Cc$, called the
tensor product, an object $I\in \Cc$ called the unit object, and
natural isomorphisms $a:\ \ot\circ (\ot\times \Cc)\to \ot\circ
(\Cc\times \ot)$ (the associativity constraint), $l:\ \ot\circ
(I\times \Cc)\to \Cc$ (the left unit constraint) and $r:\ \ot\circ
(\Cc\times I)\to \Cc$ (the right unit constraint). $a$, $l$ and
$r$ have to satisfy certain coherence conditions, we refer to
\cite[XI.2]{K} for a detailed discussion. $\Cc$ is called strict
if $a$, $l$ and $r$ are the identities on $\Cc$. McLane's
coherence Theorem asserts that every monoidal category is monoidal
equivalent to a strict one, see \cite[XI.5]{K}. The categories
that we will consider are - technically spoken - not strict, but
they can be
threated as if they were strict.\\
Let $\tau:\ \Cc\times \Cc\to\Cc\times \Cc$ be the flip functor. A prebraiding on $\Cc$
is a natural transformation $c:\ \ot\to \ot\circ \tau$ satisfying the following equations,
for all $U,V,W\in \Cc$:
$$c_{U,V\ot W}=(V\ot c_{U,W})\circ (c_{U,V}\ot W)~~;~~
c_{U\ot V,W}=(c_{U,W}\ot V)\circ (U\ot c_{V,W}).$$
$c$ is called a braiding if it is a natural isomorphism.
$c$ is called a symmetry if $c_{U,V}^{-1}=c_{V,U}$, for all $U,V\in \Cc$. We refer to
\cite[XIII.1]{K}, \cite{JS} for more details.\\
There is a natural way to associate a (pre)braided monoidal category to a monoidal
category. The weak right center $\Ww_r(\Cc)$ of a monoidal category $\Cc$ is the
category whose objects are couples of the form $(V,c_{-,V})$,
with $V\in \Cc$ and $c_{-,V}:\ -\ot V\to V\ot -$ a natural transformation
such that $c_{-, I}$ is the natural isomorphism
and
$c_{X\ot Y, V}=(c_{X, V}\ot Y)\circ (X\ot c_{Y, V})$, for all $X,Y\in \Cc$. The morphisms
are defined in the obvious way. $\Ww_r(\Cc)$ is a prebraided monoidal category;
the unit is $(I,id)$, and the tensor product is
$$
(V,c_{-, V})\ot (V',c_{-, V'})=(V\ot V', c_{-, V\ot V'})
$$
where
$$
c_{X, V\ot V'} = (V\ot c_{X, V'}) \circ (c_{X, V}\ot V')
$$
for all $X \in \Cc$. The prebraiding is given by
$$c_{V,V'}:\ (V, c_{-, V})\ot (V', c_{-, V'})\to
(V', c_{-, V'})\ot (V, c_{-, V})
$$
for all $V$, $V'\in \Cc$. The right center $\Zz_r(\Cc)$ is the
full subcategory of $\Ww_r(\Cc)$ consisting of objects
$(V,c_{-,V})$ with $c_{-,V}$ a natural isomorphism; $\Zz_r(\Cc)$
is a braided monoidal category. For more detail, we refer to
\cite[XIII.4]{K}.

\subsection{Descent data}\selabel{0.4}
Let $A$ be a commutative $k$-algebra. $\ot$ will always mean $\ot_k$, and
$A^{(n)}$ will be a shorter notation for the $n$-fold tensor product
$A\ot\cdots\ot A$. If $V$ and $W$ are right $A$-modules, then $V\ot W$ is a right
$A^{(2)}$-module.
Consider a map $g:\ A\ot V\to V\ot A$ in
$\Mm_{A^{(2)}}$. For $a\in A$ and $v\in V$, we write -
temporarily - $g(a\ot v)=\sum_i v_i\ot a_i$. Then we have the
following three maps in $\Mm_{A^{(3)}}$
\begin{equation}\eqlabel{0.4.1}
\begin{array}{ccc}
g_1:\ A\ot A\ot V\to A\ot V\ot A&;&g_1(b\ot a\ot v)=\sum_i b\ot v_i\ot a_i;\\
g_2:\ A\ot A\ot V\to V\ot A \ot A&;&g_2(a\ot b\ot v)=\sum_i v_i\ot b\ot a_i;\\
g_3:\ A\ot V\ot A\to V\ot A \ot A&;&g_3(a\ot v\ot b)=\sum_i v_i\ot a_i\ot b.
\end{array}
\end{equation}
Let $\psi:\ V\ot A\to V$ be the right $A$-action on $V$.

\begin{proposition}\prlabel{0.4.1} {\bf \cite[Prop. II.3.1]{KO}}
Assume that $g_2=g_3\circ g_1$. Then $g$ is an isomorphism if and only if
$\psi(g(1\ot v))=v$, for all $v\in V$.
\end{proposition}

In this situation, $(V,g)$ is called a descent datum. A morphism between
two descent data $(V, g)$ and
$(V', g')$ is a right $A$-linear map $f:\ V\to V'$ such that $(f\ot
A)\circ g=g'\circ (A\ot f)$. The category of descent data is
denoted by $\ul{\rm Desc}(A/k)$. We have a pair of adjoint functors $(F,G)$ between
$\Mm_k$ and $\ul{\rm Desc}(A/k)$. For $N\in \Mm_k$, $F(N)=(N\ot A,g)$, with
$g(a\ot n\ot b)=n\ot a\ot b$. $G(V,g)=\{v\in V~|~v\ot 1=g(1\ot v)\}$. The unit and counit
of the adjunction are as follows:
$$\eta_N:\ N \to (GF)(N),~~\eta_N(n)=n\ot  1;$$
$$\varepsilon_{(V,g)}:\ (FG)(V,g)=G(V,g)\ot A\to (V,g),~~\varepsilon_{(V,g)}(v\ot a)=va.$$
The Faithfully Flat Descent Theorem can now be stated as follows: if $A$ is
faithfully flat over $k$, then $(F,G)$ is an inverse pair of equivalences. This is
essentially \cite[Th\'eor\`eme 3.3]{KO}, formulated in a categorical language.
In \cite{KO}, a series of applications of descent theory are given, and there exist
many more in the literature. Also observe that the descent theory presented in \cite{KO}
is basically the affine version of Grothendieck's descent theory \cite{Gr}.

\subsection{Noncommutative descent theory and comodules over corings}\selabel{0.5}
Descent theory can be extended to the case where $A$ are
noncommutative. This was done by Cipolla in \cite{Cipolla}. After
the revival of the theory of corings initiated in
\cite{Brzezinski02}, it was observed that the results in
\cite{Cipolla} can be nicely reformulated in terms of corings.
Recall that an $A$-coring $\Cc$ is a coalgebra in the monoidal
category of $A$-bimodules. A right $\Cc$-comodule is a right
$A$-module $M$ together with a right $A$-linear map $\rho:\ M\to
M\ot_A \Cc$ satisfying appropriate coassociativity and counit
conditions. For detail on corings and comodules, we refer to
\cite{Brzezinski02,BrzezinskiWisbauer}. An important example of an
$A$-coring is Sweedler's canonical coring $\Cc=A\ot A$.
Identifying $(A\ot A)\ot_A(A\ot A)\cong A^{(3)}$, we view the
comultiplication as a map $\Delta:\ A^{(2)}\to A^{(3)}$. It is
given by the formula $\Delta(a\ot b)=a\ot 1\ot b$. The counit
$\varepsilon$ is given by $\varepsilon(a\ot b)=ab$. For a right
$A$-module $M$, we can identify $M\ot_A(A\ot A)\cong M\ot A$. A
right $A\ot A$-comodule is then a right $A$-module $V$ together
with a right $k$-linear map $\rho:\ V\to V\ot A$, notation
$\rho(v)=v_{[0]}\ot v_{[1]}$ satisfying the relations
\begin{eqnarray}
v_{[0]}v_{[1]}&=&v;\eqlabel{3.1.1}\\
\rho(v_{[0]})\ot v_{[1]}&=& v_{[0]}\ot 1\ot v_{[1]};\eqlabel{3.1.2}\\
\rho(va)&=&v_{[0]}\ot v_{[1]}a\eqlabel{3.1.3}
\end{eqnarray}
for all $v\in V$ and $a\in A$. The category of right $A\ot
A$-comodules is denoted by $\Mm^{A\ot A}$. There is an adjunction
between $\Mm_k$ and $\Mm^{A\ot A}$. Cipolla's descent data are
nothing else then $A\ot A$-comodules, and Cipolla's version of the
Faithfully Flat Descent Theorem asserts that this is a pair of
inverse equivalences if $A$ is faithfully flat over $k$,
we refer to \cite{Toronto} for a detailed discussion.\\
First observe that this machinery works for a general extension $k\to A$ of rings, that is,
$A$ and $k$ are not necessarily commutative. In this note, however, we keep $k$ commutative.\\
If $A$ is commutative, then the categories $\ul{\rm Desc}(A/k)$ and $\Mm^{A\ot A}$
are isomorphic. $(V,g)\in \ul{\rm Desc}(A/k)$ corresponds to $(V,\rho)\in \Mm^{A\ot A}$,
with $\rho(v)=g(1\ot v)$.\\
Sometimes it is argued that this generalization is not satisfactory, since there is no counterpart
to \prref{0.4.1} in the case where $A$ is noncommutative. In this note, we will present an
appropriate generalization $\ul{\rm Desc}(A/k)$ to the noncommutative situation, with
a suitable generalized version of \prref{0.4.1}, see \prref{3.4} and \reref{3.8}.

\section{The center of the category of bimodules}\selabel{2}
Throughout, $A$ is an algebra over a commutative ring $k$.

\begin{definition}\delabel{3.0}
A right Yetter-Drinfeld $A^e$-module consists of a pair $(V,\rho)$, such that $V$ is an $A$-bimodule,
$(V,\rho)\in \Mm^{A\ot A}$ and the following compatibility conditions hold:
\begin{eqnarray}
\rho(av)&=&v_{[0]}\ot av_{[1]};\eqlabel{3.1.4}\\
a\rho(v)&=&av_{[0]}\ot v_{[1]}=v_{[0]}a\ot v_{[1]}.\eqlabel{3.1.5}
\end{eqnarray}
A morphism $(V,\rho)\to (V',\rho')$ of Yetter-Drinfeld modules is
a map $f:\ V\to V'$ that is an $A$-bimodule and $A^{(2)}$-comodule
map. The category of Yetter-Drinfeld modules will be denoted by
$\YD^{A^e}$.
\end{definition}

Take $(V,\rho)\in \YD^{A^e}$. Then
\begin{equation}\eqlabel{3.1.6}
av\equal{\equref{3.1.1}}(av)_{[0]}(av)_{[1]}\equal{\equref{3.1.4}}v_{[0]}av_{[1]},
\end{equation}
and
\begin{equation}\eqlabel{3.1.7}
v_{[1]}v_{[0]}\equal{\equref{3.1.6}}v_{[0][0]}v_{[1]}v_{[0][1]}\equal{\equref{3.1.2}}v_{[0]}v_{[1]}
\equal{\equref{3.1.1}}v.
\end{equation}

\begin{proposition}\prlabel{3.2}
The forgetful functor $U:\ \YD^{A^e}\to \Mm^{A\ot A}$ is an isomorphism of categories.
\end{proposition}

\begin{proof}
We define a functor $P:\ \Mm^{A\ot A}\to \YD^{A^e}$. For $V\in
\Mm^{A\ot A}$, let $P(V)=V$ as an $A^{(2)}$-comodule, with left
$A$-action defined by $av=v_{[0]}av_{[1]}$. Then
$$\rho(av)=\rho(v_{[0]}av_{[1]})\equal{\equref{3.1.3}}\rho(v_{[0]})av_{[1]}\equal{\equref{3.1.2}}
v_{[0]}\ot av_{[1]},$$ and \equref{3.1.4} is satisfied. The left
$A$-action is associative since
$$b(av)=(av)_{[0]}b(av)_{[1]}\equal{\equref{3.1.4}}v_{[0]}bav_{[1]}=(ba)v.$$
Finally we show that \equref{3.1.5} holds:
$$av_{[0]}\ot v_{[1]}=v_{[0][0]}av_{[0][1]}\ot v_{[1]}\equal{\equref{3.1.2}}v_{[0]}a\ot v_{[1]}.$$
This shows that $P(V)\in \YD^{A^e}$. If $f:\ V\to W$ is a morphism in $\Mm^{A\ot A}$, then it is
also a morphism $P(V)\to P(W)$ in $\YD^{A^e}$. To this end, we need to show that $f$ is
left $A$-linear:
$$f(av)=f(v_{[0]}av_{[1]})=f(v_{[0]})av_{[1]}=f(v)_{[0]}af(v)_{[1]}=af(v).$$
We used the fact that $f$ is right $A$-linear and right $A\ot A$-colinear.
Finally, it is clear that the functors $P$ and $V$ are inverses.
\end{proof}

Recall from \seref{0.1} that  $\Ww_r({}_A\Mm_A)$ is the weak right
center of the monoidal category $({}_A\Mm_A, -\ot_A - , A)$ of
$A$-bimodules.

\begin{proposition}\prlabel{3.3}
The categories $\Ww_r({}_A\Mm_A)$ and $\YD^{A^e}$ are isomorphic.
\end{proposition}

\begin{proof}
Let $(V,c_{-,V})$ be an object of $\Ww_r({}_A\Mm_A)$. For every $A$-bimodule $M$,
we have an $A$-bimodule map $c_{M,V}:\ M\ot_A V\to V\ot_A M$, which is natural in $M$.
Consider
$$g=c_{A\ot A,V}:\ A^{(2)}\ot_A V\cong A\ot V\to V\ot_A A^{(2)}\cong V\ot A,$$
and define $\rho:\ V\to V\ot A$ by $\rho(v)=g(1\ot v)=v_{[0]}\ot
v_{[1]}\in V\ot A$. $c_{-,V}$ is then completely determined by
$\rho$: for $m\in M$, define the
 $A$-bimodule map $f_m:\ A^{(2)}\to M$ by the
 formula $f_m(a\ot b)=amb$.
From the naturality of $c_{-,V}$, it follows that we have a
commutative diagram
$$\xymatrix{
A^{(2)}\ot_A V\ar[d]_{f_m\ot_A V}\ar[rr]^{g}&&V\ot_A A^{(2)}\ar[d]^{V\ot_A f_m}\\
M\ot_A V\ar[rr]^{c_{M,V}}&&V\ot_AM}$$
Evaluating the diagram at $1\ot v$, we find
\begin{equation}\eqlabel{3.3.1}
c_{M,V}(m\ot_A v)=v_{[0]}\ot_A m \, v_{[1]}.
\end{equation}
We will now show that $(V,\rho)\in \YD^{A^e}$. Using the fact that
$c_{M,V}$ is right $A$-linear, well-defined and left $A$-linear,
we find
\begin{eqnarray*}
(va)_{[0]}\ot m(va)_{[1]}=c_{M,V}(m\ot_Ava)&=&c_{M,V}(m\ot_Av)a=v_{[0]}\ot_A mv_{[1]}a;\\
v_{[0]}\ot_A mav_{[1]}=c_{M,V}(ma\ot_Av)&=&c_{M,V}(m\ot_Aav)=(av)_{[0]}\ot m(av)_{[1]};\\
v_{[0]}\ot_A amv_{[1]}=c_{M,V}(am\ot_Av)&=&ac_{M,V}(m\ot_Av)=av_{[0]}\ot_A mv_{[1]}.
\end{eqnarray*}
If we take $M=A^{(2)}$ and $m=1\ot 1$ in these formulas, we obtain
\equref{3.1.3}, \equref{3.1.4} and \equref{3.1.5}. $c_{A,V}$ is
the canonical isomorphism $A\ot_A V\to V\ot_A A$, hence $v\ot_A
1=c_{A,V}(1\ot_A v)=v_{[0]}\ot_Av_{[1]}$, and \equref{3.1.1}
follows. Finally, we have the commutative diagram
$$\xymatrix{
M\ot_AN\ot_AV\ar[rr]^{c_{M\ot_A N,V}}\ar[dr]_{M\ot_A c_{N,V}}&&
V\ot_AM\ot_A N\\
&M\ot_AV\ot_A N\ar[ur]_{c_{M,V}\ot_A N}&}$$
We evaluate the diagram at $m\ot_A n\ot_A v$:
\begin{eqnarray*}
&&\hspace*{-2cm}
v_{[0]}\ot_A m\ot_A nv_{[1]}=c_{M\ot_A N,V}(m\ot_A n\ot_A v)\\
&=&
\bigl((c_{M,V}\ot_A N)\circ(M\ot_A c_{N,V})\bigr)(m\ot_A n\ot_A v)\\
&=& (c_{M,V}\ot_A N) (m\ot_A v_{[0]}\ot_A nv_{[1]}) =
v_{[0][0]}\ot_A m v_{[0][1]}\ot_A\ot_A nv_{[1]}
\end{eqnarray*}
\equref{3.1.2} follows after we take $M=N=A^{(2)}$ and $m = n = 1\ot 1$.\\
Conversely, given $(V,\rho)\in \YD^{A^e}$, we define $c_{-,V}$ using
\equref{3.3.1}. Straightforward computations show that
$(V,c_{-,V})\in \Ww_r({}_A\Mm_A)$.
\end{proof}

\begin{remark}\relabel{3.3a}
It is well-known that $A^e=A\ot A^{\rm op}$ is an $A$-bialgebroid.
The arguments in \prref{3.3} can be generalized, leading to a description of the (weak) center of
the category of modules over a bialgebroid, and to the definition of Yetter-Drinfeld module over
a bialgebroid, see \cite{sch1}. In fact, the Yetter-Drinfeld modules of \deref{3.0} are precisely the Yetter-Drinfeld
modules over the bialgebroid $A^e$, justifying our terminology.
\end{remark}

Our next aim is to show that  condition \equref{3.1.1} in \deref{3.0} can be replaced by the condition
that $g$ is invertible.

\begin{proposition}\prlabel{3.4}
Let $A$ be a $k$-algebra, and assume that $\rho:\ V\to V\ot A$ satisfies (\ref{eq:3.1.2}-\ref{eq:3.1.5}).
Then \equref{3.1.1} holds if and only if $g:\ A\ot V\to V\ot A$, $g(a\ot v)=av_{[0]}\ot v_{[1]}$ is invertible.
\end{proposition}

\begin{proof}
Assume that \equref{3.1.1} holds. For all $a\in A$ and $v\in V$, we have
\begin{eqnarray*}
&&\hspace*{-2cm}
(\tau\circ g\circ \tau \circ g)(a\ot v)=(\tau\circ g)(v_{[1]}\ot av_{[0]})\\
&=& \tau\bigl( v_{[1]}(av_{[0]})_{[0]} \ot (av_{[0]})_{[1]}\bigr)
\equal{\equref{3.1.4}}
\tau\bigl(v_{[1]}v_{[0][0]}\ot a v_{[0][1]}\bigr)\\
&\equal{\equref{3.1.5}}&
\tau\bigl(v_{[0][0]}v_{[1]}\ot a v_{[0][1]}\bigr)
\equal{\equref{3.1.2}}
\tau(v_{[0]}v_{[1]}\ot a)
\equal{\equref{3.1.1}}
a\ot v.
\end{eqnarray*}
We conclude that $\tau\circ g\circ \tau \circ g = {\rm Id}_{A\ot
V}$. Composing to the left and to the right with the switch map
$\tau$, we find $g\circ \tau\circ g\circ \tau = {\rm Id}_{V\ot
A}$. Thus $g^{-1} = \tau\circ g\circ \tau$.

Conversely, assume that $g$ is invertible. For any $v\in V$ we
have:
$$g(1\ot v_{[0]}v_{[1]})=\rho(v_{[0]}v_{[1]})\equal{\equref{3.1.3}}\rho(v_{[0]})v_{[1]}
\equal{\equref{3.1.2}} v_{[0]}\ot v_{[1]}=g(1\ot v).$$
\equref{3.1.1} follows after we apply $g^{-1}$ to both sides and
multiply the two tensor factors.
\end{proof}

\begin{proposition}\prlabel{3.5}
The (right) center of the category of $A$-bimodules coincides with
its (right) weak center: $\Zz_r({}_A\Mm_A) = \Ww_r({}_A\Mm_A)$.
\end{proposition}

\begin{proof}
Take $(V,c_{-,V})\in \Ww_r({}_A\Mm_A)$. We will show that $c_{M,V}$ is invertible, for every
$A$-bimodule $M$. Let $g$ and $\rho$ be as in \prref{3.3}. We claim that
\begin{equation}\eqlabel{3.5.1}
c_{M,V}^{-1}(v\ot_A m)=v_{[1]}m\ot_A v_{[0]}.
\end{equation}
Indeed, for all $m\in M$ and $v\in V$, we have that
$$(c_{M,V}^{-1}\circ c_{M,V})(m\ot_A v) \equal{(\ref{eq:3.3.1},\ref{eq:3.5.1})}
v_{[0][1]}m  v_{[0][0]}\ot_A v_{[0]}
\equal{\equref{3.1.2}} m\ot_A
v_{[1]}v_{[0]}\equal{\equref{3.1.7}} m\ot_A v;$$
$$(c_{M,V}\circ c_{M,V}^{-1})(v\ot_A m)\equal{(\ref{eq:3.5.1},\ref{eq:3.3.1})}
v_{[0][0]}\ot_A v_{[1]}m v_{[0][1]}
\equal{\equref{3.1.2}} v_{[0]}v_{[1]}\ot_A
m\equal{\equref{3.1.1}} v\ot_A m.$$
\end{proof}

If $V$ and $W$ are $A$-bimodules, then $V\ot W$ is an
$A^{(2)}$-bimodule. Consider a map $g:\ A\ot V\to V\ot A$ in
${}_{A^{(2)}}\Mm_{A^{(2)}}$. The maps $g_1,g_2,g_3$
defined by \equref{0.4.1} are in ${}_{A^{(3)}}\Mm_{A^{(3)}}$.

\begin{definition}\delabel{3.6}
Let $A$ be a $k$-algebra. A descent datum consists of an
$A$-bimodule $V$ together with an $A^{(2)}$-bimodule map $g:\ A\ot
V\to V\ot A$ such that $g_2=g_3\circ g_1$ and $(\psi\circ g)(a\ot
v)=v$, for all $v\in V$, where $\psi$ is the map $V\ot A\to A$,
$\psi(v\ot a)=va$. A morphism between two descent data $(V, g)$ and
$(V', g')$ is an $A$-bimodule map $f:\ V\to V'$ such that $(f\ot
A)\circ g=g'\circ (A\ot f)$. The category of descent data is
denoted by $\ul{\rm Desc}(A/k)$.
\end{definition}

\begin{proposition}\prlabel{3.7}
The categories $\ul{\rm Desc}(A/k)$ and $\YD^{A^e}$ are isomorphic.
\end{proposition}

\begin{proof}
Let $(V,\rho)\in \YD^{A^e}$, and define $g:\ A\ot V\to V\ot A$ by
$g(a\ot v)=av_{[0]}\ot v_{[1]}$. First we show that $g$ is an
$A^{(2)}$-bimodule map.
\begin{eqnarray*}
g(ba\ot cv)&=&ba(cv)_{[0]}\ot (cv)_{[1]}\equal{\equref{3.1.4}} bav_{[0]}\ot cv_{[1]}=(b\ot c)g(a\ot v);\\
g(ab\ot vc)&=&ab(vc)_{[0]}\ot (vc)_{[1]}\equal{\equref{3.1.3}}abv_{[0]}\ot v_{[1]}c\equal{\equref{3.1.5}}
av_{[0]}b\ot v_{[1]}c=g(a\ot v)(b\ot c).
\end{eqnarray*}
Now $g_3\circ g_1=g_2$ since
\begin{eqnarray*}
&&\hspace*{-2cm}
(g_3\circ g_1)(a\ot b\ot v)=g_3(a\ot bv_{[0]}\ot v_{[1]})
=a(bv_{[0]})_{[0]}\ot (bv_{[0]})_{[1]} \ot v_{[1]}\\
&\equal{\equref{3.1.4}}&av_{[0][0]}\ot bv_{[0][1]}\ot v_{[1]}
\equal{\equref{3.1.2}}av_{[0]}\ot b\ot v_{[1]}=
g_2(a\ot b\ot v).
\end{eqnarray*}
Finally, $(m\circ g)(1\ot v)=v_{[0]}v_{[1]}=v$, and
we conclude that $(V,g)\in \ul{\rm Desc}(A/k)$.\\
Conversely, let $(V,g)\in \ul{\rm Desc}(A/k)$, and define $\rho:\
V\to V\ot A$ by $\rho(v)=g(1\ot v)$. Then $f(a\ot
v)=a\rho(v)=av_{[0]}\ot v_{[1]}$. It is easy to show that
\equref{3.1.1} and (\ref{eq:3.1.3}-\ref{eq:3.1.5}) are satisfied:
\begin{eqnarray*}
v&=&(m\circ g)(1\ot v)=m(\rho(v))=v_{[0]}v_{[1]};\\
\rho(va)&=&g(1\ot va)=g(1\ot v)(1\ot a)=v_{[0]}\ot v_{[1]}a;\\
\rho(av)&=&g(1\ot av)=(1\ot a)g(1\ot v)=v_{[0]}\ot av_{[1]};\\
a\rho(v)&=&(a\ot 1)g(v)=g(a\ot v)=g(1\ot v)(a\ot 1)=v_{[0]}a\ot v_{[1]}.
\end{eqnarray*}
We have already computed $g_3\circ g_1$ and $g_2$. This
computation stays valid, since we only used \equref{3.1.4}, which
holds. Expressing that $(g_3\circ g_1)(1\ot 1\ot v)=g_2(1\ot 1\ot v)$, we
find \equref{3.1.2}. We conclude that $(V,\rho)\in \YD^{A^e}$.
\end{proof}

\begin{remarks}\relabel{3.8}
1. It follows from the proof of \prref{3.7} that the definition of
a descent datum can be restated as follows: $V\in {}_A\Mm_A$, an
invertible map $g:\ A\ot V\to V\ot A$ in
${}_{A^{(2)}}\Mm_{A^{(2)}}$ satisfying $g_2=g_3\circ g_1$.

2. We look at the particular case where $A$ is commutative. Take
$(V,g)\in \ul{\rm Desc}(A/k)$ and let $(V,\rho)$ be the
corresponding object of $\YD^{A^e}$. Then we know that
$av\equal{\equref{3.1.6}}v_{[0]}av_{[1]}=v_{[0]}v_{[1]}a\equal{\equref{3.1.1}}va$,
hence the left $A$-action on $V$ coincides with the right $A$-action. Consequently, the left and
right $A^{(2)}$-actions on $A\ot V$ and $V\ot A$ coincide. So we can view a descent datum $(V,g)$
as a right $A$-module $V$ together with a right $A^{(2)}$-linear map $g:\ A\ot V\to V\ot A$
satisfying $g_3=g_3\circ g_1$ and $(\psi\circ g)(1\ot v)=v$, or, equivalently, $g$ invertible.
These are precisely the descent data  \cite{KO} that we discussed in \seref{0.4}.
\end{remarks}

The main results of this paper are summarized as follows:

\begin{theorem}\thlabel{maincentru}
For a $k$-algebra $A$, the categories $\ul{\rm
Desc}(A/k)$, $\YD^{A^e}$, $\Mm^{A\ot A}$, $\Ww_r({}_A\Mm_A)$ and
$\Zz_r({}_A\Mm_A)$ are isomorphic.
If $A$ is faithfully flat over $k$ then these isomorphic categories
are equivalent to the category of $k$-modules.\footnote{The fact
that $\Zz_r({}_A\Mm_A)$ is equivalent to the category of
$k$-modules if $A$ is faithfully flat can be also derived from
\cite[Theorem 3.3]{sch}.}
\end{theorem}

$\Zz_r({}_A\Mm_A)$ is a braided monoidal category, hence we can
define braided monoidal structures on the five isomorphic
categories in \thref{maincentru}. In particular, the category of
comodules over the Sweedler canonical $A$-coring $A\ot A$ is
braided monoidal. Explicitly we have:

\begin{corollary}\colabel{braidedcomod}
Let $A$ be a $k$-algebra. Then $(\Mm^{A\ot A}, \, - \ot_A -, \,
A)$ is a braided monoidal category as follows: for $V\in \Mm^{A\ot
A}$, we have a left $A$-action on $V$ defined by $a\cdot v =
v_{[0]}av_{[1]}$. The tensor product is then just the tensor
product over $A$, and the coaction on $V\ot_A V'$ is given by the
formula $\rho(v\ot_A v')=v_{[0]}\ot_A v'_{[0]}\ot
v_{[1]}v'_{[1]}$. The unit is $A$, with $A\ot A$-coaction
$\rho(a)=1\ot a$. The left $A$-action on $A$ then coincides with
the left regular representation: $b\cdot a=a_{[0]}ba_{[1]}=ba$.
The braiding $c$ on $\Mm^{A\ot A}$ is given by
$$c_{V',V}(v'\ot_A v)=v_{[0]}\ot_A v'v_{[1]}~~;~~
c^{-1}_{V',V}(v\ot_A v')=v_{[1]}v'\ot_A v_{[0]}.$$
\end{corollary}

\begin{proof} This follows of course from the general theory of the center
construction, but all axioms can be easily verified directly.
\end{proof}

\begin{remark}\relabel{tbremark}
An interesting interpretation of \thref{maincentru} and
\coref{braidedcomod} was communicated to us by T. Brzezinski. In
\cite{tb06}, it was observed that there is a close relationship
between corings with a grouplike element and noncommutative
differential geometry. One of the results in this direction is the
following: the category $\Mm^{A\ot A}$ is isomorphic to the
category ${\bf Conn}(A/k,Omega(A\ot A/k)$ of right $A$-modules
with a flat connection, see \cite[Theorem 5]{tb06} or \cite[Sec.
29]{BrzezinskiWisbauer}. It then follows from \coref{braidedcomod}
that ${\bf Conn}(A/k,Omega(A\ot A/k)$ is a braided monoidal
category. In the forthcoming \cite{tb11}, the braiding
 on $\Mm^{A\ot A}$ is applied to
prove that any flat connection in a right $A$-module is an
$A$-bimodule connection.
\end{remark}

\section{Finitely generated projective algebras}\selabel{3}
Now we focus attention to the case where $A$ is finitely generated
and projective as a $k$-module, which means that the $k$-linear
map
\begin{equation}\eqlabel{3.11.1}
\varphi:\ A^* \ot A \to \Aa  =\End_k(A),~~\varphi(a^*\ot
b)(x)=\lan a^*,x\ran b
\end{equation}
is an isomorphism. Then $\varphi^{-1}({\rm Id}_A) = \sum_i
a_i^*\ot a_i$ is called a finite dual basis of $A$, and is
characterized by the formula $\sum_i \lan a_i^*,x\ran a_i=x$, for
all $x\in A$. In this situation, we also have that
\begin{equation}\eqlabel{3.11.2}
\varphi^{-1}(f)=\sum_i a_i^*\ot f(a_i),
\end{equation}
for all $f\in \Aa$. Indeed, $\varphi(\sum_i a_i^*\ot f(a_i))(x)=\sum_i \lan a_i^*,x\ran f(a_i)=f(x)$, for
all $x\in A$. Recall that we also have an algebra map
$F:\ A\ot A^{\rm op}\to \End_k(A)$, $F(a\ot b)(x)=axb$. It is then easy to show that
\begin{equation}\eqlabel{3.11.3}
\varphi(a^*\ot a)=F(a\ot 1)\circ \varphi(a^*\ot 1)=F(1\ot a)\circ \varphi(a^*\ot 1).
\end{equation}
The categories $\Mm^{A\ot A}$ and ${}_\Aa\Mm$ are isomorphic. If $V$ is a right $A\ot A$-comodule,
then we have a left $\Aa$-action given by
\begin{equation}\eqlabel{3.11.4}
f\cdot v=v_{[0]}f(v_{[1]}).
\end{equation}
for all $f \in \Aa = \End_k(A)$ and $v\in V$. Conversely, for
$V\in {}_\Aa\Mm$, we have a right $A\ot A$-coaction now given by
\begin{equation}\eqlabel{3.11.5}
\rho(v)=\sum_i f_i\cdot v\ot a_i,
\end{equation}
where we write $f_i=\varphi(a_i^*\ot 1)$. This is well-known and
can be verified easily. It also has an explanation in terms of corings: the left
dual of the $A$-coring $A\ot A$ is ${}_A\Hom(A\ot A,A)\cong
\End(A)^{\rm op}$ as $A$-rings, see for example
\cite{BrzezinskiWisbauer}. We will now transport the braided
monoidal structure of $\Mm^{A\ot A}$ to ${}_\Aa\Mm$.

If $V\in {}_\Aa\Mm$, then $V\in {}_A\Mm_A$, by restriction of
scalars via $F$. Now we also have that $V\in \Mm^{A\ot A}\cong
\YD^{A^e}$, and this gives a second $A$-bimodule structure on $V$.
These two bimodule structures coincide:
\begin{eqnarray*}
F(1\ot a)\cdot v&\equal{\equref{3.11.4}}&v_{[0]}(F(1\ot a)(v_{[1]}))=v_{[0]}v_{[1]}a=va;\\
F(a\ot 1)\cdot v&\equal{\equref{3.11.4}}&v_{[0]}(F(a\ot 1)(v_{[1]}))=v_{[0]}av_{[1]}=av.
\end{eqnarray*}
Now take $V$, $W\in {}_{\Aa}\Mm$. Then $V\ot_A W\in \Mm^{A\ot
A}\cong {}_\Aa\Mm$. We describe the $\Aa$-action on $V\ot_A W$.
\begin{eqnarray*}
&&\hspace*{-2cm}
f\cdot (v\ot_A w)\equal{\equref{3.11.4}} v_{[0]}\ot_A w_{[0]}f(v_{[1]}w_{[1]})
\equal{\equref{3.11.5}} \sum_{i,j} f_i\cdot v\ot_A (f_j\cdot w)f(a_ia_j)\\
&=& \sum_{i,j} f_i\cdot v\ot_A \bigl(F(1\ot f(a_ia_j))\circ \varphi(a_j^*\ot 1)\bigr)\cdot w\\
&\equal{\equref{3.11.3}}&\sum_{i,j} f_i\cdot v\ot_A \varphi(a_j^*\ot f(a_ia_j))\cdot w
\equal{\equref{3.11.2}}\sum_{i} f_i\cdot v\ot_A f(a_i-)\cdot w,
\end{eqnarray*}
where $f(a-)\in \Aa$ is the map sending $x\in A$ to $f(ax)$; we have an alternative description:
\begin{eqnarray*}
&&\hspace*{-2cm}
f\cdot (v\ot_A w)=\sum_{i,j} f_i\cdot v\ot_A \bigl(F(1\ot f(a_ia_j))\circ \varphi(a_j^*\ot 1)\bigr)\cdot w\\
&\equal{\equref{3.11.3}}&\sum_{i,j} f_i\cdot v\ot_A \bigl(F( f(a_ia_j)\ot 1)\circ \varphi(a_j^*\ot 1)\bigr)\cdot w\\
&=&
\sum_{i,j}\bigl(F(1\ot f(a_ia_j))\circ \varphi(a_i^*\ot 1)\bigr)\cdot v\ot_A f_j\cdot w\\
&\equal{\equref{3.11.3}}&\varphi(a_i^*\ot f(a_ia_j))\cdot v\ot_A f_j\cdot w
\equal{\equref{3.11.2}} \sum_j f(-a_j)\cdot v\ot_A f_j\cdot w.
\end{eqnarray*}
The braiding is given by the formula
$c_{V,W}(v\ot_A w)=w_{[0]}\ot_A vw_{[1]}= \sum_i f_i\cdot w\ot va_i$.
We summarize our results:

\begin{proposition}\prlabel{3.11}
Let $A$ be a finitely generated projective $k$-algebra, with
finite dual basis $\sum_i a_i^*\ot a_i$, and write
$f_i=\varphi(a_i^*\ot 1)$. The category of left
$\End_k(A)$-modules is a braided monoidal category. The tensor
product is the tensor product over $A$; a left $\End_k(A)$-module
is an $A$-bimodule by restriction of scalars via $F$. The left
$\End_k(A)$-action on $V\ot_A W$ is given by
$$
f\cdot (v\ot_A w)=\sum_{i} f_i\cdot v\ot_A f(a_i-)\cdot w=\sum_j f(-a_j)\cdot v\ot_A f_j\cdot w
$$
for all $f \in \End_k(A)$, $v \in V$ and $w\in W$. The unit object
is $A$, with its obvious left $\End_k(A)$-action $f\cdot a=f(a)$.
The braiding is given by $c_{V,W}(v\ot_A w)=\sum_i f_i\cdot w\ot_A v
\, a_i$.
\end{proposition}

\begin{remark}
As we mentioned in the introduction, the category of
Yetter-Drinfeld modules over a finite Hopf algebra is isomorphic
to the category of modules over the Drinfeld double. We have an
analogous result here: if $A$ is finite (that is, finitely
generated projective), then the category of Yetter-Drinfeld
$A^e$-modules is isomorphic to the category of representations of
$\End_k(A)$. In fact, this tells us that we can consider
$\End_k(A)$ as the Drinfeld double of $A^e$.
\end{remark}

\begin{example}\exlabel{3.12}
Let $A=k^n=\oplus_{i=1}^n ke_i$, with multiplication $e_ie_j=\delta_{ij}e_i$ and unit
$1=\sum_{i=1}^n e_i$. Let $e_i^*\in A^*$ be given by $\lan e_i^*,e_j\ran =\delta_{ij}$.
We can then identify $M_n(k)$ and $\End_k(A)$, where an endomorphism of $A$
corresponds to its matrix with respect to the basis $\{e_1,\cdots,e_n\}$. It is then easy to
see that $\varphi(e_i^*\ot e_j)=e_{ji}$, the elementary matrix with $1$ in the $(i,j)$-position
and $0$ elsewhere. Now we easily compute that
$
f_l=\varphi\bigl(\sum_r e_l^*\ot e_r\bigr)=\sum_r e_{rl}$,
$e_{ii}=F(e_i\ot 1)=F(1\ot e_{i})$ and
$e_{ij}(e_l-)= \delta_{jl}e_{ij}$.
Let $V$ and $W$ be left $M_n(k)$-modules. Then $V\ot_{k^n} W$ is again a left $M_n(k)$-module,
the left $M_n(k)$-action is given by the formulas in \prref{3.11}, which simplify as follows:
\begin{eqnarray*}
&&\hspace*{-2cm}
e_{ij}\cdot (v\ot_{k^n} w)=\sum_{l,r} e_{rl}\cdot v \ot_{k^n} \delta_{jl}e_{ij}\cdot w
= \sum_r e_{rj}\cdot v \ot_{k^n} e_{ij}\cdot w\\
&=& \sum_r e_{rj}\cdot v \ot_{k^n} (e_{ii}e_{ij})\cdot w
= \sum_r e_{rj}\cdot v \ot_{k^n} e_i(e_{ij}\cdot w)\\
&=& \sum_r (e_{rj}\cdot v)e_i \ot_{k^n} e_{ij}\cdot w
= \sum_r (e_{ii}e_{rj})\cdot v \ot_{k^n} e_{ij}\cdot w\\
&=& e_{ij}\cdot v \ot_{k^n} e_{ij}\cdot w.
\end{eqnarray*}
Finally, we compute the braiding
\begin{eqnarray*}
&&\hspace*{-2cm}
c_{V,W}(v\ot_{k^n} w)= \sum_i f_i\cdot w\ot_{k^n} ve_i
= \sum_{i,r} e_{ri}\cdot w \ot_{k^n} e_iv\\
&=& \sum_{i,r} (e_{ri}\cdot w)e_i \ot_{k^n} v
= \sum_{i,r} (e_{ii}e_{ri})\cdot w \ot_{k^n} v\\
&=& \sum_i e_{ii}\cdot w\ot_{k^n} v=w\ot_{k^n} v.
\end{eqnarray*}
The fact that the representation category of a matrix algebra is monoidal can also be understood in
a completely different way. Weak bialgebras and Hopf algebras were introduced in \cite{BNS}.
The representation category of a weak bialgebra is monoidal, see \cite{Sz,Nill,BCJ}. The tensor is
the tensor product over $H_t=\im \varepsilon_t$, where $\varepsilon_t:\ H\to H$ is given by
the formula $\varepsilon_t(h)=\lan \varepsilon,1_{(1)}h\ran 1_{(2)}$. $H=M^n(k)$ is a weak Hopf algebra,
with comultiplication and counit given by the formulas $\Delta(e_{ij})=e_{ij}\ot e_{ij}$ and
$\varepsilon(e_{ij})=1$. In fact it is a groupoid algebra, over the groupoid with $n$ objects,
and precisely one morphism $e_{ij}$ between the objects $i$ and $j$. In this situation, it is easy to
show that $\Delta(1)=\sum_l \Delta(e_{ll})=\sum_l e_{ll}\ot e_{ll}$, and $\varepsilon_t(e_{ij})=
\sum_l \lan \varepsilon, e_{ll}e_{ij}\ran e_{ll}= e_{ii}$, so that $H_t=\oplus_i ke_{ii}\cong k^n$.
The monoidal structure on $M_n(k)$ then coincides with the one that we found above. The braiding
comes from a quasitriangular structure on $M_n(k)$.
\end{example}

\section{Application to the quantum Yang-Baxter equation}\selabel{4}
Our results lead to the construction of a new family of solutions of the
quantum Yang-Baxter equation. More precisely, to every object of $\YD^{A^e}
\cong\Mm^{A\ot A}$, we can associate a solution of the quantum Yang-Baxter equation.

\begin{theorem}\thlabel{3.21}
Let $A$ be a $k$-algebra and $(V, \rho)\in \YD^{A^e}$. Then
\begin{equation}\eqlabel{3.22a}
\Omega = \Omega_{(V, \rho)} : V\ot V \to V\ot V, \qquad \Omega (v
\ot w) = w_{[0]} \ot w_{[1]} v,
\end{equation}
is a solution of the quantum Yang-Baxter
equation $\Omega^{12} \, \Omega^{13} \, \Omega^{23} = \Omega^{23}
\, \Omega^{13} \, \Omega^{12}$ in $\End (V^{(3)})$.
In particular, if $(V, \rho)\in \Mm^{A\ot A}$, then
\begin{equation}\eqlabel{3.22b}
\Omega = \Omega_{(V, \rho)} : V\ot V \to V\ot V, \qquad \Omega (v
\ot w) = w_{[0]} \ot v_{[0]} w_{[1]} v_{[1]},
\end{equation}
is a solution of the quantum Yang-Baxter
equation and $\Omega^3 = \Omega$ in $\End
(V^{(2)})$.
\end{theorem}

\begin{proof} For all $v,w,t\in V$, we have that:
\begin{eqnarray*}
&&\hspace*{-2cm}
\Omega^{12} \, \Omega^{13} \, \Omega^{23} (v \ot w \ot t) =
\Omega^{12} \,
\Omega^{13}(v \ot t_{[0]} \ot t_{[1]} \, w)\\
&{=}& \Omega^{12} \bigl((t_{[1]}w)_{[0]} \ot t_{[0]} \ot (t_{[1]}w)_{[1]} \, v \bigl)
= t_{[0][0]} \ot t_{[0][1]} \, (t_{[1]}w)_{[0]} \ot
(t_{[1]}w)_{[1]} \, v\\
&\equal{\equref{3.1.4}}& t_{[0][0]} \ot t_{[0][1]} \, t_{[1]} \,
w_{[0]} \ot w_{[1]} \, v \equal{\equref{3.1.2}} t_{[0]} \ot
t_{[1]} \, w_{[0]} \ot w_{[1]} \,v;\\
&&\hspace*{-2cm}
\Omega^{23} \, \Omega^{13} \, \Omega^{12} (v \ot
w \ot t) = \Omega^{23} \, \Omega^{13}
(w_{[0]} \ot w_{[1]} \, v \ot t)\\
&{=}& \Omega^{23}(t_{[0]} \ot w_{[1]} \, v \ot t_{[1]} \, w_{[0]})
= t_{[0]} \ot (t_{[1]}w_{[0]})_{[0]} \ot
(t_{[1]}w_{[0]})_{[1]} \, w_{[1]} \, v\\
&\equal{\equref{3.1.4}}& t_{[0]} \ot t_{[1]} \, w_{[0][0]} \ot
w_{[0][1]} \, w_{[1]} \,v \equal{\equref{3.1.2}} t_{[0]} \ot
t_{[1]} \, w_{[0]} \ot w_{[1]} \,v
\end{eqnarray*}
Thus $\Omega^{12} \, \Omega^{13} \, \Omega^{23} = \Omega^{23} \,
\Omega^{13} \, \Omega^{12}$. We have seen in \prref{3.2} that
$(V, \rho)\in \Mm^{A\ot A}$, with left $A$-action $a\cdot v =
v_{[0]} a v_{[1]}$, is an object of $\YD^{A^e}$.
With this
identification the canonical map \equref{3.22a} takes precisely
the form \equref{3.22b}. Now, for all $v,w \in V$ we
have:
\begin{eqnarray*}
\Omega^2 (v \ot w) &=& \Omega (w_{[0]} \ot v_{[0]} w_{[1]}
v_{[1]}) \equal{\equref{3.1.3}}  v_{[0][0]} \ot w_{[0][0]}
v_{[0][1]}
w_{[1]} v_{[1]} w_{[0][1]}\\
&\equal{\equref{3.1.2}} & v_{[0]} \ot w_{[0]} w_{[1]} v_{[1]}
\equal{\equref{3.1.1}}  v_{[0]} \ot w \, v_{[1]};\\
\Omega^3 (v \ot w) &=& \Omega (v_{[0]} \ot w \, v_{[1]})
\equal{\equref{3.1.3}}  w_{[0]} \ot v_{[0][0]} w_{[1]}
v_{[1]} v_{[0][1]}\\
&\equal{\equref{3.1.2}} & w_{[0]} \ot v_{[0]} w_{[1]} v_{[1]}
 = \Omega (v \ot w).
\end{eqnarray*}
\end{proof}

It is well-known that $x=x^1\ot x^2\in A^{(2)}$ is grouplike if and only if
$x^1 x^2 = 1$ and $X^1\ot X^2x^1\ot x^2=X^1\ot 1\ot X^2$. Grouplike elements
of a coring $\Cc$ are in bijective correspondence to right $\Cc$-coactions on
$A$. In the case where $\Cc=A^{(2)}$, the right $A\ot A$-coaction on $A$
associated to $x$ is $\rho(a)=x^1\ot x^2a$. For $M\in \Mm^{A\ot A}$, we define
$$M^{{\rm co}x}=\{m\in M~|~\rho(m)=mx^1\ot x^2\}.$$
$A^{{\rm co}x}$ is a subalgebra of $A$. Suppose that we have an algebra morphism
$i:\ B\to A^{{\rm co}x}$. Then we have a pair of adjoint functors (see \cite[Sec. 1]{Toronto})
$(-\ot_B A,(-)^{{\rm co}x})$ between the categories $\Mm_B$ and $\Mm^{A\ot A}$.
The right coaction on $N\ot_B A$ is simply $\rho(n\ot_B a)=n\ot_B x^1\ot x^2a$. This construction
allows us to give examples of $A\ot A$-comodules, and, a fortiori, solutions of the
quantum Yang-Baxter equation, applying \thref{3.21}. We then obtain the following.

\begin{proposition}\prlabel{3.30}
Let $x$ be a grouplike element of $A\ot A$, and let $i:\ B\to A^{{\rm co}x}$ be an
algebra morphism. For $N\in \Mm_B$, the map $\Omega:\ (N\ot_BA)^{(2)}\to (N\ot_BA)^{(2)}$
given by
\begin{equation}\eqlabel{3.30.1}
\Omega((n\ot_B a)\ot (m\ot_Bb))=(m\ot_B x^1)\ot (n\ot_B X^1x^2bX^2a)
\end{equation}
is a solution of the quantum Yang-Baxter equation.
\end{proposition}

As a particular example, we can take $x=1\ot 1$, $B=k$, $N\in \Mm_k$. Then \equref{3.30.1}
takes the form
$\Omega (m
\ot a \ot n \ot b) = n \ot 1 \ot m \ot b \, a$.
In particular, if we take $N=k$, then $\Omega(a\ot b)=
1\ot ba$, and this shows that $\Omega$ is not necessarily bijective.\\

We now present another way to construct comodules over $A\ot A$. It is shown in
\cite{ACM1} that there is a bijective correspondence between braidings on the
category of $A$-bimodules and elements $R\in A^{(3)}$
satisfying the conditions
\begin{eqnarray}\eqlabel{1.1.4}
R^1\ot aR^2\ot R^3&=&R^1\ot R^2\ot R^3a\\
\eqlabel{2.1.2}
 R^1R^2\ot R^3&=&R^2\ot R^3R^1=1\ot 1,
\end{eqnarray}
 see \cite[Theorem 2.4]{ACM1}. We then say that $(A^e,R)$ is quasitriangular, and
 we call $R$ an $R$-matrix. $R$ satisfies several other equations,
 we mention that $R$ is invariant under cyclic permutation of the tensor factors,
 and
 \begin{equation}
 R^1\ot R^2\ot 1\ot R^3=r^1R^1\ot r^2\ot r^3R^2\ot R^3\eqlabel{1.1.5},
 \end{equation}
see \cite[Theorem 2.4]{ACM1}. Yetter-Drinfeld modules can be
constructed from bimodules over quasitriangular algebras as
follows.

\begin{proposition}\prlabel{nouAna}
Let $A$ be a $k$-algebra, let $V$ be an $A$-bimodule, and let $R\in A^{(3)} $ be an $R$-matrix.
Consider $\rho_R:\ V\to V\ot A$, $\rho_R(v) = R^{1} \, v \, R^{2} \ot R^{3} = v_{[0]} \ot v_{[1]}$.
Then $(V, \rho_R)\in \YD^{A^e}$, and the associated solution of the quantum Yang-Baxter
equation is
$$
\Omega_R = \Omega_{(V, \rho_R)} : V\ot V \to V\ot V, \quad
\Omega_R (v \ot w) = R^1 w R^2 \ot R^3 v.
$$
\end{proposition}

\begin{proof}
We show that $(V,\rho_R)\in \Mm^{A\ot A}$, that is, $\rho$ satisfies
(\ref{eq:3.1.1}-\ref{eq:3.1.2}).
 \equref{3.1.1} follows from \equref{2.1.2}.
\equref{3.1.2} is equivalent to
$$(R^{1} \, v \, R^{2})_{[0]} \ot (R^{1} \, v \, R^{2})_{[1]} \ot R^{3} = R^{1} \, v \, R^{2} \ot 1 \ot R^{3}$$
and to
\begin{equation}\eqlabel{ana1}
r^{1} \, R^{1} \, v \, R^{2} \, r^{2} \ot r^{3} \ot R^{3} = R^{1}
\, v \, R^{2} \ot 1 \ot R^{3}.
\end{equation}
 Using (\ref{eq:1.1.4}-\ref{eq:1.1.5}), we obtain:
$$
R^{1} \ot R^{2} \ot 1 \ot R^{3} \equal{\equref{1.1.5}} r^{1} \,
R^{1} \ot r^{2} \ot r^{3} \, R^{2} \ot R^{3}
\equal{\equref{1.1.4}} r^{1} \, R^{1} \ot R^{2} \, r^{2} \ot r^{3}
\ot R^{3}
$$
and \equref{ana1} follows. It follows from \prref{3.2} that $(V,\rho)\in\YD^{A^e}$,
and we are done if we can show that the left $A$-action on $V$ given by \equref{3.1.6}
coincides with the original left $A$-action. This can be shown easily:
$$v_{0]}av_{[1]}=R^1vR^2aR^3=aR^1vR^2R^3=av.$$
We used (\ref{eq:1.1.4}-\ref{eq:2.1.2}), combined with the fact that $R$ is invariant under
cyclic permuation of the tensor factors.
\end{proof}

Several examples of $R$-matrices are presented in \cite{ACM1}. In particular,
if $A$ is an Azumaya algebra, then we have a unique $R$-matrix.
Applying \prref{nouAna} to \cite[Example 2.8]{ACM1}, we obtain the following.

\begin{example}
Let $A = M_n(k)$ be a matrix algebra and $V$ an $M_n(k)$-bimodule.
Then the map
$$
\Omega : V \ot V \to V \ot V, \quad \Omega (v \ot w) =
 \sum_{i,j,k=1}^n e_{ij} \, w \, e_{ki} \ot e_{jk} \, v
$$
is a solution of the quantum Yang-Baxter equation. $e_{ij}$
is the elementary matrix with $1$ in the $(i,j)$-position and $0$
elsewhere.
\end{example}

\begin{problem}
Let $V$ be a finite dimensional vector
space over a field $k$ and let $\Omega \in \End (V^{(2)})$ be a solution
of the quantum Yang-Baxter equation such that $\Omega^3 =
\Omega$. Does there exist an algebra $A$ and a right $A\ot A$-coaction
on $V$ such that $\Omega$ is given by \equref{3.22b}?
\end{problem}

\section*{Acknowledgment} We would like to thank Gabriella B\"ohm and Tomasz
Brzezi\'nski for their comments on the first version of this paper.

\end{document}